\documentclass{amsart}
\usepackage{amssymb}
\usepackage{ifthen}
\usepackage{graphicx}
\usepackage{xcolor}

\newtheorem{thm}{Theorem}
\newtheorem{cor}{Corollary}
\newtheorem{lem}{Lemma}

\newtheorem{rem}{Remark}

\newtheorem{conj}{Conjecture}
\theoremstyle{definition}
\newtheorem{example}[equation]{Example}
\newtheorem{prob}[equation]{Problem}

\newcommand{\IC}{{\mathbb C}}
\newcommand{\ID}{{\mathbb D}}

\newcommand{\D}{{\mathbb D}}

\def\be{\begin{equation}}
\def\ee{\end{equation}}

\newcommand{\bee}{\begin{enumerate}}
\newcommand{\eee}{\end{enumerate}}

\newcommand{\blem}{\begin{lem}}
\newcommand{\elem}{\end{lem}}
\newcommand{\bthm}{\begin{thm}}
\newcommand{\ethm}{\end{thm}}
\newcommand{\bcor}{\begin{cor}}
\newcommand{\ecor}{\end{cor}}
\newcommand{\beg}{\begin{example}}
\newcommand{\eeg}{\end{example}}
\newcommand{\begs}{\begin{examples}}
\newcommand{\eegs}{\end{examples}}
\newcommand{\bdefe}{\begin{defin}}
\newcommand{\edefe}{\end{defin}}
\newcommand{\bprob}{\begin{prob}}
\newcommand{\eprob}{\end{prob}}
\newcommand{\bei}{\begin{itemize}}
\newcommand{\eei}{\end{itemize}}

\newcommand{\bcon}{\begin{conj}}
\newcommand{\econ}{\end{conj}}
\newcommand{\bcons}{\begin{conjs}}
\newcommand{\econs}{\end{conjs}}
\newcommand{\bprop}{\begin{propo}}
\newcommand{\eprop}{\end{propo}}
\newcommand{\br}{\begin{rem}}
\newcommand{\er}{\end{rem}}
\newcommand{\brs}{\begin{rems}}
\newcommand{\ers}{\end{rems}}
\newcommand{\bo}{\begin{obser}}
\newcommand{\eo}{\end{obser}}
\newcommand{\bos}{\begin{obsers}}
\newcommand{\eos}{\end{obsers}}
\newcommand{\bpf}{\begin{pf}}
\newcommand{\epf}{\end{pf}}
\newcommand{\ba}{\begin{array}}
\newcommand{\ea}{\end{array}}
\newcommand{\beq}{\begin{eqnarray}}
\newcommand{\beqq}{\begin{eqnarray*}}
\newcommand{\eeq}{\end{eqnarray}}
\newcommand{\eeqq}{\end{eqnarray*}}

\begin{document}

\title[Univalency of certain transform of univalent functions]{Univalency of certain transform of univalent functions}

\author[M. Obradovi\'{c}]{Milutin Obradovi\'{c}}
\address{Department of Mathematics,
Faculty of Civil Engineering, University of Belgrade,
Bulevar Kralja Aleksandra 73, 11000, Belgrade, Serbia}
\email{obrad@grf.bg.ac.rs}

\author[N. Tuneski]{Nikola Tuneski}
\address{Department of Mathematics and Informatics, Faculty of Mechanical Engineering, Ss. Cyril and
Methodius University in Skopje, Karpo\v{s} II b.b., 1000 Skopje, Republic of North Macedonia.}
\email{nikola.tuneski@mf.edu.mk}

\subjclass[2000]{ 30C45}
\keywords{analytic, univalent, transform}

\maketitle

\begin{abstract}
We consider univalency problem in the unit disc $\ID$ of the function
 $$g(z)=\frac{(z/f(z))-1}{-a_{2}},$$
where $f$ belongs to some classes of univalent functions in $\ID$  and   $a_{2}=\frac{f''(0)}{2}\neq 0$.
\end{abstract}

\medskip

\section{Introduction}

Let ${\mathcal A}$ denote the family of all analytic functions $f$
in the unit disk $\ID := \{ z\in \IC:\, |z| < 1 \}$  satisfying the normalization
$f(0)=0= f'(0)-1$, i.e., $f$ has the form
\be\label{eq 1}
f(z)=z+a_{2}z^{2}+a_{3}z^{3}+\ldots.
\ee

 Let $\mathcal{S}$, $\mathcal{S}\subset \mathcal{A}$, denote the class of univalent functions in $\ID$, let $\mathcal{S}^{\star}$ be the subclass of ${\mathcal A}$ (and $\mathcal{S}$ which are starlike in $\ID$ and let $\mathcal{U} $ denote the set of all  $f\in {\mathcal A}$  satisfying
the condition
\begin{equation}\label{exp-U}
\left |\operatorname{U}_f(z) \right | < 1 \qquad (z\in\ID),
\end{equation}
where
\be\label{eq 5}
\operatorname{U}_{f}(z):=\left(\frac{z}{f(z)}\right)^{2}f'(z)-1.
\ee

\medskip

In \cite[Theorem 4]{OT-2019} the authors consider the problem  of univalency for the function
\be\label{eq 2}
g(z)=\frac{(z/f(z))-1}{-a_{2}},
\ee
where $f\in\mathcal{U}$  has the form \eqref{eq 1} with $a_{2}\neq 0$.
They proved the following

\medskip

\noindent\bf{Theorem A}. \it
Let $f\in\mathcal{U}$. Then, for the function $g$ defined by the expression \eqref{eq 2} we have
\begin{itemize}
  \item[(a)] $|g'(z)-1|<1$ for $|z|<|a_{2}|/2 $;
  \item[(b)] $g\in \mathcal{S}^{\star}$  in the disk $|z|<|a_{2}|/2$, and even more
$$ \left|\frac{zg'(z)}{g(z)}-1 \right|<1 $$
in the same disk ;
\item[(c)] $g\in \mathcal{U}$ in the disk $|z|<|a_{2}|/2$ if $0<|a_{2}|\leq 1.$
\end{itemize}
These results are the best possible. \rm

\medskip

For the proof of the previous theorem the authors used the next representation for the class  $\mathcal{U}$
(see \cite{OPR-1996} and \cite{OP-2011}). Namely,  if $f\in\mathcal{U}$ then
\be\label{eq 3}
\frac{z}{f(z)} = 1-a_2 z -z\omega(z),
\ee
where function $\omega$ is analytic in $\D$ with  $|\omega(z)|\leq |z|<1$ for all $z\in\D$.
The appropriate function $g$ from \eqref{eq 2} has the form
\be\label{eq 4}
g(z)=z+\frac{1}{a_{2}}z\omega(z).
\ee

\medskip

\section{Results}

In this paper we consider other cases of Theorem A(c) and certain related results.

\medskip

\bthm \label{th 1}
Let $f\in\mathcal{U}$. Then the function $g$ defined by the equation \eqref{eq 2}
belongs to $\mathcal{U}$ in the disc
$$ |z|<\sqrt{\frac{1-|a_2|+\sqrt{|a_{2}|^{2}+2|a_{2}|-3}}{2}},$$
i.e., satisfies \eqref{exp-U} on this disc, if $ \frac{5}{4}\leq|a_{2}|\leq 2$.
\ethm
\begin{proof} For the first part of the proof we use the same method as in \cite{OT-2019}.
By the definition of the class $\mathcal{U},$ i.e., inequality \eqref{exp-U}, and using the next estimation for the function $\omega $,
$$|z\omega'(z)-\omega(z)|\leq\frac{r^{2}-|\omega(z)|^{2}}{1-r^{2}},$$ (where $|z|=r$ and
$|\omega(z)|\leq r$),
after some calculations,  we obtain
\[
\begin{split}
|\operatorname{U}_{g}(z)|&=
\left|\frac{\frac{1}{a_{2}}\left[z\omega'(z)-\omega(z)\right]-\frac{1}{a^{2}_{2}}\omega^{2}(z)}{\left[1+\frac{1}{a_{2}}\omega_{1}(z)\right]^{2}}\right|\\
&\leq\frac{|a_{2}|\cdot|z\omega'(z)-\omega(z)|+|\omega(z)|^{2}}{\left(|a_{2}|-|\omega(z)|\right)^{2}}\\
&\leq \frac{|a_{2}|\cdot \frac{r^{2}-|\omega(z)|^{2}}{1-r^{2}}+|\omega(z)|^{2}}{\left(|a_{2}|-|\omega(z)|\right)^{2}}\\
&=:\frac{1}{1-r^{2}}\cdot \varphi(t).
\end{split}
\]
Here,
\be\label{eq 6}
\varphi(t)=\frac{|a_{2}|r^{2}-(|a_{2}|-1+r^{2})t^{2}}{(|a_{2}|-t)^{2}}
\ee
and $|\omega(z)|=t$, $0\leq t\leq r$. From here we have that
$$\varphi'(t)=\frac{2|a_{2}|}{(|a_{2}|-t)^{3}}\cdot \left[r^{2}-(|a_{2}|-1+r^{2})t\right]$$
(where $|a_{2}|-t>0$ since $|a_{2}|\geq \frac{5}{4}>1>t$).
Next, $\varphi'(t)=0$ for
$$t_{0}=\frac{r^{2}}{|a_{2}|-1+r^{2}}$$
and $0 \leq t_{0}\leq r $ if
$$\frac{r^{2}}{|a_{2}|-1+r^{2}}\leq r,$$
which is equivalent to
$$r^{2}-r+|a_{2}|-1 \geq 0. $$
The last relation is valid for $ \frac{5}{4}\leq|a_{2}|\leq 2$ and every $0\leq t<1$. It means that the maximal value of the function $\varphi$ on $[0,r]$ is
$$ \varphi(t_{0})=\frac{(|a_{2}|-1+r^{2})r^{2}}{(|a_{2}|-1)(|a_{2}|+r^{2})}.$$
Finally,
$$|\operatorname{U}_{g}(z)|\leq \frac{1}{1-r^{2}}\cdot \varphi(t_{0})=
\frac{(|a_{2}|-1+r^{2})r^{2}}{(1-r^{2})(|a_{2}|-1)(|a_{2}|+r^{2})}<1$$
if
$$ r^{4}-(1-|a_{2}|)r^2+(1-|a_{2}|)< 0,$$
or if
$$ r<\sqrt{\frac{1-|a_2|+\sqrt{|a_{2}|^{2}+2|a_{2}|-3}}{2}}.$$
This completes the proof.
\end{proof}

\medskip

For our next consideration we need the next lemma.
\blem \label{le_1}
Let $f\in \mathcal{A}$ be of the form \eqref{eq 1}. If
\be\label{eq 7}
\sum_{2}^{\infty}n|a_{n}|\leq1,
\ee
then
\[
\begin{split}
 |f'(z)-1|&<1 \qquad (z\in\ID),\\[2mm]
\left|\frac{zf'(z)}{f(z)}-1\right|&<1\qquad (z\in\ID)
\end{split}
\]
(i.e. $f\in \mathcal{S}^{\star}$ ), and
$f\in \mathcal{U}$.
\elem

For the proof of $f\in \mathcal{U}$ in the lemma see \cite{OP-2011}, while the rest easily follows.

\medskip

Further, let $\mathcal{S}^{+}$ denote the class of univalent functions in the unit disc
with the representation
\be\label{eq 8}
\frac{z}{f(z)}=1+b_{1}z+b_{2}z^{2}+\ldots, \quad b_{n}\geq 0,\,\,n=1,2,3,\ldots.
\ee
For example, the Silverman class (the class with negative coefficients) is included
in the class $\mathcal{S}^{+}$, as well as the Koebe function $k(z)=\frac{z}{(1+z)^{2}}\in\mathcal{S}^{+}$.
The next characterization is valid for the class $\mathcal{S}^{+}$ (for details see \cite{obpo-2009})
\be\label{eq 9}
f\in \mathcal{S}^{+} \quad \Leftrightarrow \quad \sum _{n=2}^{\infty}(n-1)b_{n}\leq 1.
\ee

\medskip

\bthm \label{th 2}
Let $f\in\mathcal{S}^{+}$. Then the function $g$ defined by \eqref{eq 2} belongs to the class
$\mathcal{U}$ in the disc $|z|<|a_{2}|/2$ and the result is the best possible.
\ethm

\begin{proof}
Using the representation \eqref{eq 8}, the corresponding function $g$ has the form
$$g(z)=\frac{\frac{z}{f(z)}-1}{-a_{2}}=\frac{\frac{z}{f(z)}-1}{b_{1}}
=z+\sum_{2}^{\infty} \frac{b_{n}}{b_{1}}z^{n}\quad (b_{1}\neq 0), $$
and from here
$$\frac{1}{r}g(rz)=z+\sum_{2}^{\infty}\frac{b_{n}}{b_{1}}r^{n-1}z^{n}\quad (0<r\leq 1). $$
Then, after applying Lemma \ref{le_1}, we have
\[
\begin{split}
\sum_{2}^{\infty}n|a_{n}|&=\sum_{2}^{\infty} n\frac{b_{n}}{b_{1}}r^{n-1}\\
&=\frac{1}{b_{1}}\sum_{2}^{\infty} (n-1)b_{n}\frac{n}{n-1}r^{n-1}\\
&\leq \frac{2r}{b_{1}}\sum_{2}^{\infty}(n-1)b_{n}\leq \frac{2r}{b_{1}}\leq 1
\end{split}
\]
if $r\leq \frac{b_{1}}{2}=\frac{|a_{2}|}{2}$. It means,  by the same lemma, that
$g\in\mathcal{U}$ in the disc $|z|<|a_{2}|/2.$ \\

In order to show that the result is the best possible, let consider the function $f_{1}$ defined by
\be\label{eq 10}
\frac{z}{f_{1}(z)}=1+bz+z^{2},\quad 0< b\leq2.
\ee
Then, $f_{1}\in\mathcal{S}^{+}$ is of type $f_1(z)=z-bz^2+\cdots$,  so the function  $$ g_{1}(z)=\frac{\frac{z}{f_{1}(z)}-1}{b}=z+\frac{1}{b}z^{2}$$ is such that
$$\left|\left(\frac{z}{g_{1}(z)}\right)^{2}g'_{1}(z)-1\right|\leq\frac{\frac{1}{b^{2}}|z|^{2}}{\left(1-\frac{1}{b}|z|\right)^{2}}<1$$
when  $|z|<b/2$. This implies that $g_{1}$
belongs to the class $\mathcal{U}$ in the disc  $|z|<b/2$.
On the other hand, since
$g'_{1}(-b/2)=0$, the function $g_{1}$ is not univalent in a bigger disc, implying that the result is the best possible.
\end{proof}

\medskip

\bthm \label{th 3}
Let $f\in\mathcal{S}$. Then the function $g$ defined by \eqref{eq 2} belongs to the class
$\mathcal{U}$ in the disc $|z|<r_{0}$,  where $r_{0}$ is the unique real root of the equation
\be\label{eq 11}
\frac{3r^{2}-2r^{4}}{(1-r^{2})^{2}}-\ln(1-r^{2})=|a_{2}|^{2}
\ee
on the interval $(0,1)$.
\ethm
\begin{proof}
We apply the same method as in the proof of the previous theorem. Namely, if $f\in\mathcal{S}$ has the representation \eqref{eq 8}, then
\be\label{eq 12}
\sum _{n=2}^{\infty}(n-1)|b_{n}|^{2}\leq 1
\ee
(see \cite{Go}, Theorem 11, p. 193, Vol. 2). Also, using  \eqref{eq 2},  \eqref{eq 8} and \eqref{eq 12}, we have $a_2=-b_1$, and
$$\frac{1}{r}g(rz)=z+\sum_{2}^{\infty}\frac{b_{n}}{b_{1}}r^{n-1}z^{n}, \quad 0<r\leq 1. $$
So,
\[
\begin{split}
\sum_{n=2}^{\infty}n|a_{n}|&= \sum_{n=2}^{\infty} n\frac{|b_{n}|}{|b_{1}|}r^{n-1}\\
&= \frac{1}{|b_{1}|}\sum_{n=2}^{\infty} \sqrt{n-1}\cdot |b_{n}|\cdot \frac{n}{\sqrt{n-1}}\cdot r^{n-1}\\
&\leq \frac{1}{|b_{1}|}\cdot \left(\sum_{n=2}^{\infty}(n-1)|b_{n}|^{2}\right)^{1/2}
\cdot \left(\sum_{n=2}^{\infty}\frac{n^{2}}{n-1}r^{2(n-1)}\right)^{1/2}\\
&\leq \frac{1}{|b_{1}|}\left(r^{2}\sum_{n=2}^{\infty}(n-1)(r^{2})^{n-2}+2r^{2}\sum_{n=2}^{\infty}(r^{2})^{n-2}
+\sum_{n=2}^{\infty}\frac{1}{n-1}(r^{2})^{n-1}\right)^{1/2} \\
&= \frac{1}{|b_{1}|}\left[\frac{3r^{2}-2r^{4}}{(1-r^{2})^{2}}-\ln(1-r^{2})\right]^{1/2}\leq 1
\end{split}
\]
if  $|z|<r_{0}$, where $r_{0}$ is the root of the equation
$$\frac{3r^{2}-2r^{4}}{(1-r^{2})^{2}}-\ln(1-r^{2})=|b_{1}|^{2} (=|a_{2}|^{2}). $$
We note that the function on the left side of this equation is an increasing one on the interval $(0,1)$, so the equation has a unique root when $0<|a_{2}|\leq 2.$
\end{proof}

\end{document}